\documentclass[12pt,twoside]{article}
\usepackage{graphicx}
\usepackage{amsfonts}
\usepackage{amsbsy}
\usepackage{amssymb}
\usepackage{amsmath,tikz,pgfplots}
\usetikzlibrary{decorations.pathmorphing}
\usepackage{cite}
\usepackage{graphics}
\def\bpsp{\begin{pspicture}}
\def\epsp{\end{pspicture}}

\newtheorem{theorem}{Theorem}[section]
\newtheorem{remark}[theorem]{Remark}
\newtheorem{example}[theorem]{Example}
\newtheorem{lemma}[theorem]{Lemma}
\newtheorem{corollary}[theorem]{Corollary}
\newtheorem{definition}[theorem]{Definition}
\newtheorem{proposition}[theorem]{Proposition}
\newtheorem{note}{Note}
\newtheorem{case}{Case}
\newtheorem{conjecture}{Conjecture}
\newtheorem{question}{Question}

\newcommand{\bea}{\begin{eqnarray}}
\newcommand{\eea}{\end{eqnarray}}
\newcommand{\beq}{\begin{eqnarray*}}
\newcommand{\eeq}{\end{eqnarray*}}

\def\m4{\mbox{\rm ~(mod $4$)}}

\def \bd{\begin{definition}}
\def \ed{\end{definition}}
\def \bqu{\begin{question}}
\def \equ{\end{question}}
\def \bcc{\begin{conjecture}}
\def \ecc{\end{conjecture}}
\def \bt{\begin{theorem}}
\def \et{\end{theorem}}
\def \bl{\begin{lemma}}
\def \el{\end{lemma}}
\def \bc{\begin{corollary}}
\def \ec{\end{corollary}}
\def \be{\begin{equation}}
\def \ee{\end{equation}}
\def \ben{\begin{enumerate}}
\def \een{\end{enumerate}}
\def \ba{\begin{array}}
\def \ea{\end{array}}
\def \bp{\begin{proposition}}
\def \ep{\end{proposition}}
\def \bx{\begin{example}}
\def \ex{\end{example}}
\def \br{\begin{remark}}
\def \er{\end{remark}}
\def \bdsc{\begin{description}}
\def \edsc{\end{description}}

\def \bn{\begin{case}}
\def \en{\end{case}}
\def \bnt{\begin{note}}
\def \ent{\end{note}}
\def\1{1\!\!1}

\def\mm2{\mbox{\rm ~(mod $2$)}}
\def\m4{\mbox{\rm ~(mod $4$)}}

\def\qed{\nolinebreak\hfill\rule{.2cm}{.2cm}\par\addvspace{.5cm}}

\def\m{\mu}

\def\1{\textbf{1}}
\def\0{\textbf{0}}

\begin{document}
\title{On  adjacency and Laplacian  cospectral non-isomorphic signed graphs }
\author{ Tahir Shamsher$^{a}$, S. Pirzada$^{b}$, Mushtaq A. Bhat$^{c}$,\\
$^{a,b}$ {\em Department of Mathematics, University of Kashmir, Srinagar, Kashmir, India}\\
$^{c}${\em  Department of Mathematics, National Institute of Technology, Srinagar, India}\\
$^a$tahir.maths.uok@gmail.com,~~ $^b$pirzadasd@kashmiruniversity.ac.in,\\$^c$mushtaqab@nitsri.net }
\date{}

\pagestyle{myheadings} \markboth{Shamsher, Pirzada, Bhat}{On  adjacency and Laplacian  cospectral non-isomorphic signed graphs }
\maketitle
\vskip 5mm
\noindent{\footnotesize \bf Abstract.} Let $\Gamma=(G,\sigma)$ be a signed graph, where $\sigma$ is the sign function on the edges of $G$. In this paper, we use the operation of partial transpose to obtain  non-isomorphic Laplacian cospectral signed graphs. We will introduce two new operations on signed graphs. These operations will establish a relationship between the adjacency spectrum of one signed graph with the Laplacian spectrum of another signed graph. As an application, these new operations  will be utilized to construct several pairs of cospectral non-isomorphic signed graphs. Finally, we construct  integral signed graphs.

\vskip 3mm



\section{Introduction}\label{sec1}

Let $G=(V(G),E(G))$ be a simple connected graph with vertex set $V(G)=\{v_1,v_2,\dots,v_n\}$ and edge set $E(G)=\{e_1,e_2,\dots,e_m\}.$ A signed graph is defined to be a pair $\Gamma=(G,\sigma)$, with $G=(V(G),E(G))$ as the underlying graph and $\sigma: E(G) \rightarrow \{-1, 1\}$ as the signing function. In this manuscript, bold lines denote positive edges, and dashed lines denote negative edges. If all the edges of a signed graph have either a positive or negative sign, then it is said to be homogenous; otherwise, it is said to be heterogeneous. Signed graphs are a generalization of graphs, since they are homogeneous signed graphs with each edge positive. The sign of a signed graph is the product of the signs of its edges. A signed graph is said to be positive (or negative) if its sign is positive (or negative), that is, it has an even (or odd) number of negative edges. If all of the edges of a signed graph are positive (or negative), it is said to be all-positive (or all-negative) respectively. A signed graph is said to be balanced if none of its cycles is negative, otherwise unbalanced.\\
\indent In a signed graph $\Gamma=(G,\sigma)$, the degree of a vertex $v$ is the same as its degree in the underlying graph $G$ (denoted by $d_v(G)$). Let $\Gamma$ be a signed graph with vertex set $V(G)$ and let $X\subset V(G)$ be a nonempty set. Let $\Gamma^X$ denote the signed graph obtained from $\Gamma$ by reversing signs of edges between $X$ and $V(G)-X$. Then, we say $\Gamma^X$ is switching equivalent to $\Gamma$. Here, we note that the switching is an equivalence relation and preserves the eigenvalues of the adjacency and the Laplacian matrix including their multiplicities.  A switching class is represented by a single signed graph.\\
\indent The adjacency matrix of a signed graph $\Gamma$, whose vertices are ${v_1,v_2,\cdots,v_n}$, is the $n\times n$ matrix $A(\Gamma)=(a_{ij})$, where
\begin{equation*}
a_{ij}=\left\{\begin{array}{lr}\sigma(v_i,v_j), &\mbox{if there is an edge from $v_i$ to $v_j,$}\\
0, &\mbox {otherwise.}
\end{array}\right.\end{equation*}
For a graph $G$, the Laplacian matrix is $L(G) = D(G) - A(G)$ and signless Laplacian matrix is $Q(G)=D(G) + A(G)$, where $A(G)$ and $D(G)$ are respectively the adjacency matrix and the diagonal matrix of vertex degrees of $G$. The Laplacian matrix of $\Gamma$ is $ L(\Gamma) = L(G, \sigma ) = D(G) - A(\Gamma)$.  Note that $L(G, +) = L(G)$ and $L(G, -) = Q(G)$.   The characteristic polynomial $|xI-A(\Gamma)|$ and eigenvalues  of the adjacency matrix $A(\Gamma)$ of the signed graph $\Gamma$ are denoted by $\phi_\Gamma(x)$ and $\lambda_1,  \lambda_2, \cdots ,\lambda_n$ respectively. The characteristic polynomial $|xI-L(\Gamma)|$ and eigenvalues  of the Laplacian matrix $L(\Gamma)$ of the signed graph $\Gamma$ are denoted by $\psi_\Gamma(x)$ and $\mu_1,  \mu_2, \cdots , \mu_n$ respectively. For a graph $G$ (resp. signed graph $\Gamma$), eigenvalues of its adjacency matrix and Laplacian matrix are called adjacency and
Laplacian eigenvalues of $G$ (resp. $\Gamma$). Clearly, $A(\Gamma)$ and $L(\Gamma)$  are real symmetric and so all their eigenvalues are real. Let the signed graph $\Gamma$ of order $n$ has distinct eigenvalues $\lambda_1,\lambda_2,\cdots,\lambda_k$ and let their respective multiplicities be $m_1,m_2,\cdots,m_k$. The adjacency spectrum of $\Gamma$ is written as $Spec(\Gamma)=\{\lambda^{(m_1)}_1,\lambda^{(m_2)}_2,\cdots,\lambda^{(m_k)}_k\}$. A signed graph is said to be an integral signed graph if its adjacency spectrum consists of integers only.\\
\indent  Given a graph G,  its subdivision graph ${\bf \mathit{ S}}(G)$ is obtained from $G$ by replacing each of its edge by a path of length $2$, or, equivalently, by inserting an additional vertex into each edge of $G$. If two signed graphs have the same adjacency spectrum (Laplacian spectrum), they are said to be cospectral (Laplacian cospectral); otherwise, they are noncospectral (Laplacian noncospectral). Any two isomorphic signed graphs are cospectral (Laplacian cospectral). A signed graph is said to be determined by its adjacency spectrum if cospectral signed graphs are isomorphic signed graphs. It is well-known that in general the adjacency spectrum does not determine the signed graph, and this problem has attracted to identify, if any, cospectral non-isomorphic signed graphs for a given class of signed graphs. For open problems in signed graphs we refer to \cite{openproblem}.\\
 The energy of a graph $G$ is the sum of the absolute values of its adjacency eigenvalues. This concept was extended to signed graphs by Germina, Hameed and Zaslavsky \cite{ghz}.  The energy of a signed graph $\Gamma$ with eigenvalues $x_1,x_2,\dots,x_n$ is defined as $\mathcal{E}(\Gamma)=\sum_{j=1}^{n}| x_j|$. Two signed graphs of same order are said to be equienergetic if they have the same  energy.\\
\indent Harary {\cite{h}} pioneered the use of signed graphs in connection with the study of social balance theory.  Signed graphs have been intensively explored in a variety of fields such as group theory, topological graph theory and classical root system.  The reader is referred to {\cite{Z}} for a complete bibliography on signed graphs.

The rest  of the paper is organized as follows. In Section $2$, we give some preliminary results which will be used in the sequel. In Section $3$, we define the concept of partial transpose to signed graphs and use it to obtain  non-isomorphic Laplacian cospectral signed graphs. In section $4$, we introduce two new operations on signed graphs and which will be utilized to construct cospectral non-isomorphic signed graphs, noncospectral equienergetic signed graphs and integral signed graphs.
\section{Preliminaries}\label{sec2}
In this section, we recall some previously established results which will be required in the subsequent sections.
\begin{definition}\label{2.1} {\em \cite{ds}} Let $P=(p_{ij})\in M_{m\times n}(\mathbb{R})$ and $Q$ $\in M_{p \times q}(\mathbb{R})$. The Kronecker product of $P$ and $Q$, denoted by $P\otimes Q$, is defined as
\begin{equation*}
P\otimes Q=
\begin{pmatrix}
   p_{11}Q &  p_{12}Q & \dots &  p_{1n}Q\\
   p_{21}Q &  p_{22}Q & \dots &  p_{2n}Q\\
   \vdots  &  \vdots & \ddots &  \vdots\\
    p_{m1}Q &  p_{m2}Q & \dots &  p_{mn}Q
\end{pmatrix}.
\end{equation*}
\end{definition}
 \begin{lemma}\label{2.2}{\em \cite{ds}} Let $P$, $Q$ $\in M_n(\mathbb{R})$ be two square matrices of order $n$. Let $\lambda$ be an
eigenvalue of matrix $P$ with corresponding eigenvector $x$ and $\mu$ be an eigenvalue of
matrix $Q$ with corresponding eigenvector $y$. Then $\lambda  \mu $ is an eigenvalue of $P \otimes Q$ with
corresponding eigenvector $x \otimes y$.
\end{lemma}
  \begin{lemma}\label{2.3}{\em \cite{ghz}} Let $ \Gamma_1 $ and $ \Gamma_2 $ be two signed graphs with respective eigenvalues $ x_1, x_2, \dots, x_{n_1}$ and $ y_1, y_2, \dots, y_{n_2}$. Then\\
   $(i)$ the eigenvalues of $ \Gamma_1 \times \Gamma_2$ are $x_i+y_j$, for all $i=1, 2,\dots, n_1$ and  $j=1, 2,\dots, n_2$,\\
$(ii)$ the eigenvalues of $\Gamma_1 \otimes \Gamma_2$ are $x_i y_j$, for all $i=1, 2,\dots, n_1$ and  $j=1, 2,\dots, n_2$.
\end{lemma}
\indent The Cartesian product (or sum) of two signed graphs $\Gamma_1=(V(G_1),E(G_1),\sigma_1)$ and  $\Gamma_2=(V(G_2), E(G_2),\sigma_2)$, denoted by $\Gamma_1\times \Gamma_2$, is the signed graph $(V(G_1)\times V(G_2),E,\sigma)$, where the edge set is that of the Cartesian product of underlying unsigned graphs and the sign function is defined by
 \begin{equation*}\sigma((u_i,v_j),(u_k,v_l))=\left \{\begin{array}{lr}\sigma_1(u_i,u_k), &\mbox{if $j=l$},\\
 \sigma_2(v_j,v_l), &\mbox{if $i=k.$}
 \end{array} \right.\end{equation*}
\indent The Kronecker product (or conjunction) of two signed graphs  $\Gamma_1=(V(G_1),E(G_1),\sigma_1)$ and  $\Gamma_2=(V(G_2),E(G_2),\sigma_2)$, denoted by $\Gamma_1\otimes \Gamma_2$, is the signed graph $(V(G_1)\otimes V(G_2),E,\sigma)$, where the edge set is that of the Kronecker product of underlying unsigned graphs and the sign function is defined by $\sigma((u_i,v_j),(u_k,v_l))=\sigma_1(u_i,u_k)\sigma_2 (v_j,v_l).$
\begin{lemma}\label{2.4}{\em \cite{mpe}}
 Let $\Gamma$ be an unbalanced signed graph with at least one edge, whose spectrum is symmetric about the origin,  having eigenvalues $\xi_1, \xi_2, \dots, \xi_n$. Then $\Gamma\times K_2$ and $\Gamma\otimes K_2$ are unbalanced, noncospectral and equienergetic if and only if $|\xi_j| \ge 1$, for all $j=1,2,\dots,n.$
\end{lemma}
\begin{lemma}\label{2.5}{\em \cite{hj}}
Let $P(\bullet)$ be a given polynomial. If $\mu$ is an eigenvalue of $A\in M_{n}$, while $y$ is an associated eigenvector, then $P(\mu)$ is an eigenvalue of the matrix $P(A)$ and $y$ is an associated eigenvector with $P(\mu)$.
\end{lemma}

 \begin{lemma}\label{2.6} {\em \cite{mpe}} Let $\Gamma$ be a signed graph of order n. Then the following statements are equivalent.\\
$(i)$ The spectrum of $\Gamma$ is symmetric about the origin,\\
$(ii)$ $\phi_\Gamma(x)=x^n+\sum_{k=1}^{\lfloor{\frac{n}{2}}\rfloor}(-1)^k b_{2k}x^{n-2k}$, where $b_{2k}$ are non negative integers for all $k=1,2,\dots, {\lfloor{\frac{n}{2}}\rfloor} $, \\
$(iii)$ $\Gamma$ and $-\Gamma$ are cospectral, where $-\Gamma$ is the signed graph obtained by negating sign of each edge of $\Gamma$.
\end{lemma}
\begin{lemma}\label{2.7} {\em \cite{r}} For infinitely many $n$, there exists a family of $2^k$ pairwise
nonisomorphic Laplacian integral, Laplacian cospectral graphs on $n$
vertices, where $k >\frac{n}{(2 log_2(n))}. $
\end{lemma}

\section{Constructing Laplacian cospectral non-isomorphic signed graphs}

 Dutta {\cite{3}} constructed large families of non-isomorphic  signless Laplacian
cospectral graphs using partial transpose on graphs. In this section, we define partial transpose for signed graphs.  Let $\Gamma=(G,\sigma)$ be a signed graph on $2n$ vertices with vertex set $V(G)=V_{1} \cup V_{2}$, such that $V_{1} \cap V_{2}=\emptyset$, and
$V_{1}=\left\{u_{1}, u_{2}, \ldots u_{n}\right\}, V_{2}=\left\{v_{1}, v_{2}, \ldots v_{n}\right\}.$
 We denote $\left\langle V_{1}\right\rangle_{\Gamma}$ and $\left\langle V_{2}\right\rangle_{\Gamma}$ as the induced signed subgraphs of $\Gamma$ generated by $V_{1}$ and $V_{2}$, respectively. The spanning signed subgraph of $\Gamma$ consisting of the signed edge set $\left\{\left(u_{i}, v_{j}\right)\in E(\Gamma): u_{i} \in V_{1}, v_{j} \in V_{2}\right\}$ is denoted by $\left\langle V_{1}, V_{2}\right\rangle_{\Gamma}$. Let $\widehat{E\left(\left\langle V_{1}, V_{2}  \right\rangle_{\Gamma}\right)}$ $=\left\{\left(u_{j}, v_{i}\right),\right.$ where $\left.\left(u_{i}, v_{j}\right) \in E\left(\left\langle V_{1}, V_{2}\right\rangle_{\Gamma}\right)\right\} $ be the set of edges, which  suggests that given any signed edge $\left(u_{i}, v_{j}\right) \in E\left(\left\langle V_{1}, V_{2}\right\rangle_{\Gamma}\right)$ there is a  unique signed edge $\left(u_{j}, v_{i}\right)\in$ $\widehat{E\left(\left\langle V_{1}, V_{2}\right\rangle_{\Gamma}\right)}$ with the same sign as the sign of edge $\left(u_{j}, v_{i}\right)$  in $\left\langle V_{1}, V_{2}\right\rangle_{\Gamma}$. It is easy to see that the set  $E\left(\left\langle V_{1}, V_{2}\right\rangle_{\Gamma}\right)$ consists of all edges of the form $\left(u_{i}, v_{j}\right)$ in $\Gamma$, but the existence of an edge $\left(u_{i}, v_{j}\right)$ in $\Gamma$ does not assure the existence of $\left(u_{j}, v_{i}\right)$ in $\Gamma$.\\
 \indent The partial transpose of a signed graph $\Gamma$, denoted by $\Gamma^{\tau}$, is defined as $\Gamma^{\tau}=\Gamma-E\left(\left\langle V_{1}, V_{2}\right\rangle_{\Gamma}\right)+\widehat{E\left(\left\langle V_{1}, V_{2}\right\rangle_{\Gamma}\right)}$. Note that, subtracting $E\left(\left\langle V_{1}, V_{2}\right\rangle_{\Gamma}\right)$ indicates to remove all the existing signed edges in $\Gamma$ of the form $\left(u_{i}, v_{j}\right) \in E\left(\left\langle V_{1}, V_{2}\right\rangle_{\Gamma}\right)$. Then we include the non-existing edges $\left(u_{j}, v_{i}\right) \in \widehat{E\left(\left\langle V_{1}, V_{2}\right\rangle_{\Gamma}\right)}$ to construct $\Gamma^{\tau}$. If $i=j,$ then the edge $\left(u_{i}, v_{i}\right)$ will be removed and added again, that is the edge $(u_{i}, v_{i})$ is unaltered under partial transpose. Therefore, partial transpose of a signed graph $\Gamma$ is an operation on the edge set which replaces the signed edge $\left(u_{i}, v_{j}\right)$ with the sign $\sigma= \pm 1$, with the corresponding signed edge $\left(u_{j}, v_{i}\right)$ with the same sign $\sigma$. \\
\indent Consider the signed graphs $\Gamma_1$ and $\Gamma_1^{\tau}$ as shown in Figure 1. Here, we have $V_1=\{u_1,u_2,u_3\}$, $V_2=\{v_1,v_2,v_3\}$ and $E\left(\left\langle V_{1}, V_{2}\right\rangle_{\Gamma_1}\right)= \{(u_1,v_1),(u_1,v_3)\}$. Thus, $\widehat{E\left(\left\langle V_{1}, V_{2}\right\rangle_{\Gamma_1}\right)}= \{(u_1,v_1),(u_3,v_1)\}$. Here, we replace the existing signed edge $(u_1, v_3)$ with the non-existing signed edge $(u_3, v_1)$.\\

\noindent{\bf Remark 3.1} The partial transpose of a signed graph is labelling dependent. Therefore, isomorphic signed graphs may have non-isomorphic partial transposes, depending on the labellings. The partial transpose keeps $\left\langle V_{1}\right\rangle$ and $\left\langle V_{2}\right\rangle$ unaltered. The total number of vertices remains unchanged. If degree of a vertex $w_{i}$ in the signed graph $\Gamma$ is $d_{\Gamma}(w_i)$, then
\begin{equation*}\sum\limits_{i=1}^{n} (d_{\Gamma}(u_i)+d_{\Gamma}(v_i))= \sum\limits_{i=1}^{n} (d_{\Gamma^\tau}(u_i)+d_{\Gamma^\tau}(v_i)).
\end{equation*}
\indent A cycle $C_l^\sigma(v_1, v_2,\cdots, v_l, v_1)$ in a signed graph $\Gamma=(G,\sigma)$ is a finite sequence of distinct vertices such that $(v_i, v_{i+1}) \in E(\Gamma)$ for all $i = 1, 2, \cdots, l-1$ and $(v_l , v_1) \in E(\Gamma)$. We denote the negative edges in the signed cycle $C_l^\sigma(v_1, v_2,\cdots, v_l, v_1)$ by putting the bar over the corresponding adjacent vertices. For example, the cycle $C_4^\sigma(v_1,v_2,v_3,v_4,v_1)$  on four vertices  such that the only edge $(v_1,v_2)\in E(\Gamma)$ has negative sign will be denoted by $C_4^-(\overline{v_1,v_2},v_3,v_4,v_1)$. Similarly if only two consecutive  edges $(v_1,v_2), (v_2,v_3)\in E(\Gamma)$ have negative signs, then the cycle $C_4^\sigma(v_1,v_2,v_3,v_4,v_1)$ will be denoted by $C_4^+(\overline{v_1,v_2,v_3},v_4,v_1)$. In a signed graph $\Gamma$, a signed $TU$-subgraph $\Gamma'$  is a signed subgraph whose components are trees or unbalanced
unicyclic graphs, namely the unique cycle contains an odd number of negative edges. Thus, if
$H$ is a signed $TU$-subgraph, then $H = T_1 \cup T_2 \cup \cdots \cup T_p \cup U_1 \cup U_2 \cup \cdots \cup U_q$, where
 $T_i’s$ are trees and  $U_i’s$ are unbalanced unicyclic graphs. The weight of the signed
$TU$-subgraph H is defined as $w(H)=4^{q} \prod_{i=1}^{p}\left|T_{i}\right|$, where $|T_{i}|$ is the number of vertices in the tree $T_i$. Note that we define $\prod_{i=1}^{p}\left|T_{i}\right|=1$ when $p=0$. The relation between the coefficients of the Laplacian characteristic polynomial with the $TU$-subgraphs of a signed graph can be seen in [\cite{bc},Theorem 3.9]. Let $\Gamma$ be a signed graph with Laplacian characteristic polynomial $\psi(\Gamma, x)=x^{n}+a_{1} x^{n-1}+\cdots+a_{n-1} x+a_{n}$. Then its coefficients are given by
\begin{equation}
a_{i}=(-1)^{i} \sum_{H \in \mathcal{H}_{i}(\Gamma)} w(H) \quad(i=1,2, \ldots, n),
\end{equation}
where $\mathcal{H}_i(\Gamma)$ denotes the set of signed $TU$-subgraphs of $\Gamma$  containing $i$ edges. Two sets of signed $TU$-subgraphs $\mathcal{H}_{i}(\Gamma)$ and $\mathcal{H}_{i}\left(\Gamma'\right)$ are comparable if
\begin{equation*}
\sum_{H \in \mathcal{H}_{i}(\Gamma)} w(H)=\sum_{H \in \mathcal{H}_{i}\left(\Gamma'\right)} w(H).
\end{equation*}
 Now, Eq. $(3.1)$ suggests that if $\Gamma$ and $\Gamma'$ are  Laplacian cospectral, then the sets of their signed $TU$-subgraphs are comparable for all $i=1,2, \ldots m$, where $m$ is the number of edges in the signed graph $\Gamma$. We say two signed graphs $\Gamma_1$ and $\Gamma_2$ are comparable if $\mathcal{H}_{i}(\Gamma_1)$ and $\mathcal{H}_{i}(\Gamma_2)$ are comparable for all $i .$ As an example, two signed paths with equal number of vertices are comparable. \\
{\bf Example 3.1} Consider the signed graphs $\Gamma_1$ and
$\Gamma_1^{\tau}$ as shown in Figure 1.
 We observe that $\Gamma_1$ contains two cycles $C_3^+(u_1,u_2,u_3,u_1)$ and $C_4^-(\overline{u_1,v_1},v_2,v_3,u_1)$. The partial transpose $\Gamma_1^\tau$  of $\Gamma_1$ which is obtained by replacing the signed edge $(u_1,v_3)$ with $(u_3,v_1)$ preserves both these cycles $C_3^+(u_1,u_2,u_3,u_1)$ and $C_4^-(\overline{u_1,v_1},u_3,u_2,u_1)$ in $\Gamma_1^\tau$. As balanced cycle do not contribute to the coefficients of the Laplacian characteristic polynomial of a signed graph, therefore the signed $TU$-subgraphs generated by  $C_4^-(\overline{u_1,v_1},v_2,v_3,u_1)$ and $C_4^-(\overline{u_1,v_1},u_3,u_2,u_1)$ in $\Gamma_1$  and $\Gamma_1^\tau$ respectively, are isomorphic, and have same contribution in $\psi_{\Gamma_1}(x)$ and $\psi_{\Gamma_1^\tau}(x)$.\\
\indent The signed edges $K_{1,3}=\{(u_1,v_1), (u_1,v_3), (u_1, u_3)\}$ form a tree in $\Gamma_1$. It is replaced by an unbalanced unicyclic $TU$-subgraph $C_3^-(\overline{u_1, v_1},u_3,u_1)$ in $\Gamma_1^\tau$. Clearly, the signed $TU$-subgraphs generated by $K_{1,3}$  and $C_3^-(\overline{u_1, v_1},u_3,u_1)$ in $\Gamma_1$  and $\Gamma_1^\tau$ respectively have equal contribution in $\psi_{\Gamma_1}(x)$ and $\psi_{\Gamma_1^\tau}(x)$. Therefore, all the signed $TU$-subgraphs of  $\Gamma_{1}$ and $\Gamma_{1}^\tau$ are comparable. Thus they have same Laplacian characteristic polynomial, which can be easily calculated by Eq. $(3.1)$ and is given by
 \begin{equation*}  \psi_{\Gamma_{1}}(x)=\psi_{\Gamma_1^\tau}(x)=x^{6}-14 x^{5}+73 x^{4}-176 x^{3}+196 x^{2}-88 x+12.
 \end{equation*}
 \begin{figure}
\centering
	\includegraphics[scale=.8]{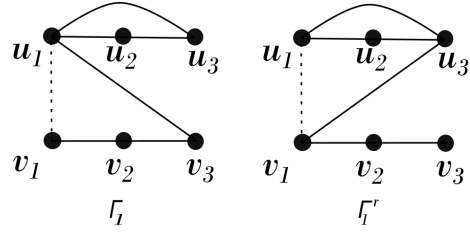}
	\caption{Signed graph $\Gamma_1$ and its partial transpose $\Gamma_1^\tau$.}
	\label{Figure 1}
\end{figure}

 If we add a signed edge $(w_i,w_j): $ $w_i,w_j \in V(G)$ in a signed graph $\Gamma=(G,\sigma)$, then the resultant signed graph will be denoted by $\Gamma'=\Gamma + \{(w_i,w_j)\}$. Similarly, $\Gamma'=\Gamma - \{(w_i,w_j)\}$ denotes the signed graph obtained by removing an edge $(w_i,w_j):$ $w_i,w_j \in V(G)$ in a signed graph $\Gamma=(G,\sigma)$. Whether the added/removed edge $(w_i,w_j)$ is positive or negative, we denote a negative edge by  $\overline{(w_i,w_j)}$, and a positive edge without a bar over the edge $(w_i,w_j)$.
 \begin{theorem}{\label{3.1}}
 Let the signed subgraphs $\left\langle V_{1}\right\rangle_{\Gamma}$ and $\left\langle V_{2}\right\rangle_{\Gamma}$ of the signed graph $\Gamma$ be two paths on $n$ vertices with each edge being positive. Let $\left\langle V_{1}, V_{2}\right\rangle_{\Gamma}$ be an empty signed graph. For the  new constructed signed graph $\Gamma_1=\Gamma +\left\{\left(u_{1}, u_{n}\right),\overline{\left(u_{1}, v_{1}\right)},\left(u_{1}, v_{n}\right)\right\}$, \\
 $(i)$ the signed graph $\Gamma_1$ is non-isomorphic and Laplacian cospectral to its partial transpose $\Gamma_1^\tau$.\\
 $(ii)$ the signed graphs $\Gamma_2=\Gamma_1-\left\{ {\left(u_{n-1}, u_{n}\right)},\overline{\left(u_{1}, v_{1}\right)}\right\}+\left\{ \overline{\left(u_{n-1}, u_{n}\right)},{\left(u_{1}, v_{1}\right)}\right\}$ and $\Gamma_3=\Gamma_1^\tau-\left\{(u_{n-1}, u_{n})\right\}+\left\{ \overline{\left(u_{n-1}, u_{n}\right)}\right\}$ are non-isomorphic and Laplacian cospectral.
 \end{theorem}
{\bf Proof.} $(i)$ The cycles generated by additional three edges and their incidence with existing edges in $\Gamma_1$ are $C_n^+\left(u_{1}, u_{2}, u_{3}, \ldots u_{n}, u_{1}\right)$ and $C_{n+1}^-\left(\overline{u_{1}, v_{1}}, v_{2}, \ldots v_{n}, u_{1}\right)$, respectively. The signed spanning subgraph $\left\langle V_{1}, V_{2}\right\rangle_{\Gamma}$ contains only two signed edges which are $\overline{(u_{1}, v_{1})}$ and $\left(u_{1}, v_{n}\right) .$ Partial transpose replaces $\left(u_{1}, v_{n}\right)$ with $\left(u_{n}, v_1\right).$ The cycles $C_n^+\left(u_{1}, u_{2}, u_{3}, \ldots u_{n}, u_{1}\right)$ and \\$C_{n+1}^-\left(\overline{u_{1}, v_{1}}, v_{2}, \ldots v_{n}, u_{1}\right)$ remain invariant (into an isomorphic cycles) under partial transpose on $\Gamma_1$. As  balanced cycles do not contribute to the coefficients of the Laplacian characteristic polynomial of a signed graph, therefore the signed $TU$-subgraphs generated by  $C_{n+1}^-\left(\overline{u_{1}, v_{1}}, v_{2}, \ldots v_{n}, u_{1}\right)$ and $C_{n+1}^-\left(\overline{u_{1}, v_{1}}, u_{n}, \ldots u_{2}, u_{1}\right)$ in $\Gamma_1$ and $\Gamma_1^\tau$, respectively,  are isomorphic and have same contribution in $\psi_{\Gamma_1}(x)$ and $\psi_{\Gamma_1^\tau}(x)$. Now,  the edges $\{\overline{(u_1,v_1)}, (u_1,v_n), (u_1, u_n)\}$ form a signed tree in $\Gamma_1$. It is replaced by an unbalanced unicyclic $TU$-subgraph $C_3^-(\overline{u_1, v_1},u_n,u_1)$ in $\Gamma_1^\tau$ and have same contribution in $\psi_{\Gamma_1}(x)$ and $\psi_{\Gamma_1^\tau}(x)$. Therefore, all the signed $TU$-subgraphs of  $\Gamma_{1}$ and $\Gamma_{1}^\tau$ are comparable. Thus, by Eq. $(3.1)$, they have same Laplacian characteristic polynomial, which proves the result in this case.\\
\begin{figure}
\centering
	\includegraphics[scale=.8]{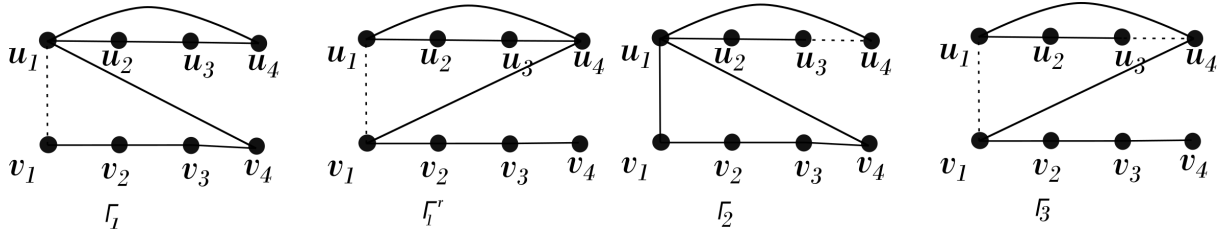}
	\caption{Signed graphs $\Gamma_1$, $\Gamma_1^\tau$, $\Gamma_2$ and $\Gamma_3$.}
	\label{Figure 2}
\end{figure}
$(ii)$ The proof of (ii) is similar to that of $(i)$.\qed

\noindent{\bf Example 3.2}$~$Consider the signed graphs $\Gamma_{1}$, $\Gamma_{1}^\tau$, $\Gamma_{2}$  and $\Gamma_{3}$ as given in Figure 2. They are constructed by using Theorem \ref{3.1}. Their Laplacian characteristic polynomials are respectively given as below.\\
\begin{equation*}\psi_{\Gamma_{1}}(x)=\psi_{\Gamma_{1}^{\tau}}(x)=x^{8}-18 x^{7}+131 x^{6}-498 x^{5}+1061 x^{4}-1256 x^3+764x^2-200x+16,\end{equation*}
\begin{equation*}\psi_{\Gamma_{2}}(x)=\psi_{\Gamma_{3}}(x)=x^{8}-18 x^{7}+131 x^{6}-498 x^{5}+1065 x^{4}-1288 x^3+848x^2-280x+36.\end{equation*}
Clearly, the signed graphs $\Gamma_{1}$ and $\Gamma_{1}^\tau$  are non-isomorphic and Laplacian cospectral. Also, $\Gamma_{2}$ and $\Gamma_{3}$  are non-isomorphic and Laplacian cospectral signed graphs.

\begin{theorem}{\label{3.2}}
 Let the signed subgraphs $\left\langle V_{1}\right\rangle_{\Gamma}$ and $\left\langle V_{2}\right\rangle_{\Gamma}$ of  $\Gamma$ be two cycles on $n$ vertices with each edge being positive. Let $\left\langle V_{1}, V_{2}\right\rangle_{\Gamma}$ be an empty signed graph. Given two non-adjacent vertices $u_{i}$ and $u_{j}$ with $i<j$,  construct a new signed graph $\Gamma_1=\Gamma +\left\{\left(u_{i}, u_{j}\right),\overline{\left(u_{i}, v_{i}\right)},\left(u_{i}, v_{j}\right)\right\}$. Then\\
 $(i)$ the signed graph $\Gamma_1$ is non-isomorphic and Laplacian cospectral to its partial transpose $\Gamma_1^{\tau}$,\\
 $(ii)$ the signed graphs $\Gamma_2=\Gamma_1-\left\{ {\left(u_{n-1}, u_{n}\right)},\overline{\left(u_{i}, v_{i}\right)},\left(u_{j-1}, u_{j}\right) \right\}+\left\{ \overline{\left(u_{n-1}, u_{n}\right)},{\left(u_{1}, v_{1}\right)}, \overline{\left(u_{j-1}, u_{j}\right)}\right\}$ and $\Gamma_3=\Gamma_1^\tau-\left\{ {\left(u_{n-1}, u_{n}\right)},\left(u_{j-1}, u_{j}\right) \right\}+\left\{ \overline{\left(u_{n-1}, u_{n}\right)}, \overline{\left(u_{j-1}, u_{j}\right)}\right\}$ are non-isomorphic and Laplacian cospectral,\\
  \begin{figure}
\centering
	\includegraphics[scale=.6]{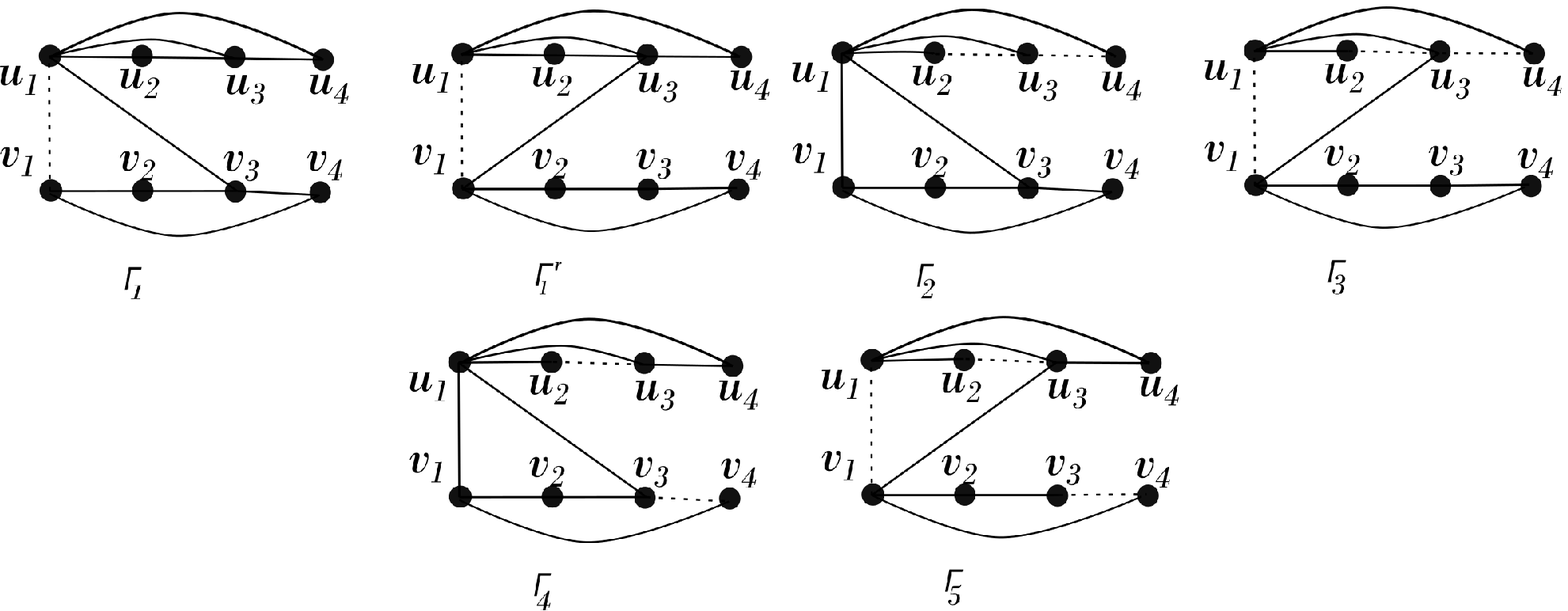}
	\caption{Signed graphs $\Gamma_1$ , $\Gamma_1^\tau$, $\Gamma_2$, $\Gamma_3$, $\Gamma_4$ and $\Gamma_5$.}
	\label{Figure 3}
\end{figure}
 $(iii)$ the signed graphs $\Gamma_4=\Gamma_1-\left\{ {\left(v_{n-1}, v_{n}\right)},\overline{\left(u_{i}, v_{i}\right)},\left(u_{j-1}, u_{j}\right) \right\}+\left\{ \overline{\left(v_{n-1}, v_{n}\right)},{\left(u_{1}, v_{1}\right)}, \overline{\left(u_{j-1}, u_{j}\right)}\right\}$ and $\Gamma_5 =\Gamma_1^\tau-\left\{ {\left(v_{n-1}, v_{n}\right)},\left(u_{j-1}, u_{j}\right) \right\}+\left\{ \overline{\left(v_{n-1}, v_{n}\right)}, \overline{\left(u_{j-1}, u_{j}\right)}\right\}$ are non-isomorphic and Laplacian cospectral.
  \end{theorem}
 {\bf Proof.} The proof is similar to that of Theorem \ref{3.1}. \qed

\noindent{\bf Example 3.3}\label{3.3}  Consider the  signed graphs $\Gamma_1$, $\Gamma_1^\tau$, $\Gamma_2$, $\Gamma_3$, $\Gamma_4$ and $\Gamma_5$ as shown in Figure 3. Here $n=4$, $i=1$ and $j=3$. The signed graph $\Gamma_1$, which is generated by Theorem \ref{3.2} is non-isomorphic and Laplacian cospectral to its partial transpose $\Gamma_1^\tau$. We obtain the signed graph $\Gamma_2$ from $\Gamma_1$ by replacing the positive edges $(u_2,u_3)$ and  $(u_3,u_4)$ with negative edges $\overline{(u_2,u_3)}$ and  $\overline{(u_3,u_4)}$ and negative edge $\overline{(u_1,v_1)}$ with the positive edge $(u_1,v_1)$. Also, the signed graph $\Gamma_3$ is obtained from $\Gamma_1^\tau$ by replacing the positive edges $(u_2,u_3)$ and  $(u_3,u_4)$ with negative edges $\overline{(u_2,u_3)}$ and  $\overline{(u_3,u_4)}$. The signed graphs $\Gamma_2$ and $\Gamma_3$ are non-isomorphic and Laplacian cospectral.  Similarly the non-isomorphic and Laplacian cospectral signed graphs $\Gamma_4$ and $\Gamma_5$ are obtained from $\Gamma_1$ and $\Gamma_1^\tau$, respectively, as  in  Theorem \ref{3.2}.\\
{\bf Remark 3.2} In Example $3.3$, we have seen that $\Gamma_1$ and $\Gamma_1^\tau$ are Laplacian cospectral signed graphs. Also, we have mentioned that $\Gamma_1^\tau$ is the partial transpose of $\Gamma_1$. But, not all signed graphs are Laplacian cospectral to their partial transpose, for instance, consider the signed graphs $\Gamma$ and $\Gamma^\tau$ as given in Figure $4$. It is easy to calculate that the Laplacian characteristic polynomials of $\Gamma$ and $\Gamma^\tau$ are
 $ \psi_{\Gamma}(x) = x^{6}-12 x^{5}+52 x^{4}-105 x^{3}+104 x^{2}-48 x+8$ and $ \psi_{\Gamma^{\tau}}(x) = x^{6}-12 x^{5}+51 x^{4}-94 x^{3}+72 x^{2}-18 x.$ \\
\indent Let $G$ be a graph and $\Gamma=(G, \sigma)$ be a signed graph on $G$. Hou et al. \cite{y} raised the following two problems. \\
{\bf Problem 1.} Let $G$ be a graph, $\Gamma_{1}=\left(G, \sigma_{1}\right), \Gamma_{2}=\left(G, \sigma_{2}\right)$ be two signed graphs on $G$, and $\operatorname{det} L\left(\Gamma_{1}\right)=\operatorname{det} L\left(\Gamma_{2}\right).$ Are $L\left(\Gamma_{1}\right)$ and $L\left(\Gamma_{2}\right)$ cospectral?\\
{\bf Problem 2.} Do there exist pairs $\Gamma_{1}=\left(G_{1}, \sigma_{1}\right), \Gamma_{2}=\left(G_{2}, \sigma_{2}\right)$ of signed graphs that have either of the following properties (i) and (ii)?\\
(i) $\Gamma_{1}$ and $\Gamma_{2}$ are not balanced but Laplacian cospectral such that $G_{1}$ and $G_{2}$ are nonisomorphic.\\
(ii) $\Gamma_{1}$ and $\Gamma_{2}$ are not balanced but Laplacian cospectral such that $G_{1}$ and $G_{2}$ are
not cospectral.\\
\indent The statement of Problem 1 is not always true. To see this, let $\Gamma_1$ be a signed graph as shown in Figure 3.  Let $\Gamma'$ be the signed graph obtained from $\Gamma_1$  by replacing the negative edge $\overline{(u_1,v_1)}$ with positive edge $(u_1,v_1)$ and positive edge $(v_3,v_4)$ with negative edge $\overline{(v_3,v_4)}.$ The Laplacian characteristic polynomials of $\Gamma_1$  and $\Gamma'$  are respectively given by\\
$~~~~~~~~~\psi_{\Gamma_{1}}(x)=x^{8}-22 x^{7}+197 x^{6}-928 x^{5}+2476 x^{4}-3736 x^3+2976x^2-1056x+128,$\\
$~~~~~~~~~ \psi_{\Gamma'}(x)=x^{8}-22 x^{7}+197 x^{6}-928 x^{5}+2476 x^{4}-3748 x^3+3048x^2-1152x+128.$\\
 The underlying graphs of  $\Gamma_1$  and $\Gamma'$ are isomorphic and  $\operatorname{det} L\left(\Gamma_{1}\right)=\operatorname{det} L\left(\Gamma'\right).$  It is clear that the signed graphs $\Gamma_1$  and $\Gamma'$ are not Laplacian cospectral and this answers Problem 1.\\
\indent For Problem 2, consider the signed graph $\Gamma_1$  and its partial transpose  $\Gamma_1^{\tau}$ as given in Figure 3. Clearly, the underlying graphs of  $\Gamma_1$  and $\Gamma_1^{\tau}$ are non-isomorphic. The  unbalanced signed graphs $\Gamma_1$  and   $\Gamma_1^{\tau}$ are Laplacian cospectral. Also, it is easy to see that the underlying graph of $\Gamma_1$  and   $\Gamma_1^{\tau}$ are not cospectral and this answers Problem 2.
  \begin{figure}
\centering
	\includegraphics[scale=.6]{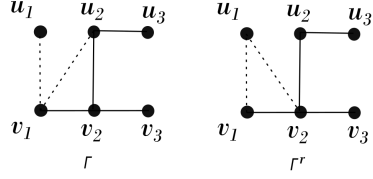}
	\caption{Signed graphs $\Gamma$ and $\Gamma^\tau$.}
	\label{Figure 4}
\end{figure}

\section{Constructing  cospectral(adjacency) non-isomorphic signed graphs, integral signed graphs and equienergetic(adjacency) signed graphs}\label{sec3}

The novel cospectral non-isomorphic signed graph constructions have implications for the complexity of the graph isomorphism problem. This necessitates the creation of methods for detecting and/or creating cospectral non-isomorphic graphs. Seidel switching, Godsil–McKay (GM) switching, and others are well-known approaches for constructing cospectral graphs. In 2019, Belardo et al. {\cite{bcs}} used the Godsil-Mckay-type procedures developed for graphs to construct the pairs of cospectral switching non-isomorphic signed graphs. In this section, we will introduce two new operations in signed graphs. These operations establish the relationship of the adjacency spectrum of one signed graph with the Laplacian spectrum of another signed graph. Furthermore, these operations will be utilized to construct the pairs of cospectral non-isomorphic signed graphs and integral signed graphs. \\
\indent The usual orientation of edges in digraphs differs slightly from the orientation of signed graphs. In fact in signed graphs, instead of one arrow, we can use two arrows assigned to edges. Bidirected graphs are the result of this. An orientated signed graph, more exactly, is an ordered pair $\Gamma_{\vartheta}=(\Gamma,\vartheta)$, where
\begin{equation}
\vartheta:V(G)\times E(G)\rightarrow \{0,1,-1\}
\end{equation}
satisfying the following three conditions.\\
$(a)$ $\vartheta(u, vw)=0$ whenever $u \neq v,w;u,v,w \in V(G)$ and $vw\in E(G)$,\\
$(b)$ $\vartheta(v, vw)=1$ ( or $-1$) if an arrow at $v$ is going into (rep. out of) $v$ . For illustration, see Figure 5,\\
$(c)$ $\vartheta(v, vw)$ $\vartheta(w, vw)= -\sigma(vw).$ \\

\begin{figure}
\centering
	\includegraphics[scale=.6]{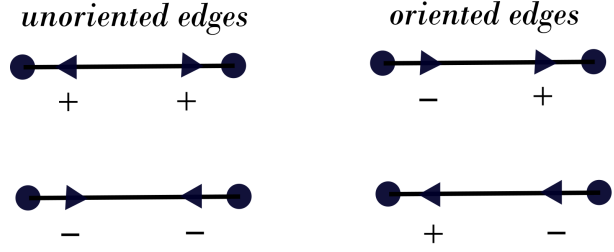}
	\caption{Bidirected edges in signed graphs.}
	\label{Figure 5}
\end{figure}

As a result, positive edges are oriented edges, whereas negative edges are unoriented (see  Figure 5). Therefore, every bidirected graph is also a signed graph. The converse is likewise true, however, one arrow (at any end) can be taken at random, whereas the other arrow (in light of $(c)$ above) cannot.
 For an oriented signed graph $\Gamma _{\vartheta}$, its incidence matrix $B_{\vartheta}=(b_{ij})$ is a matrix, whose rows correspond to vertices and columns to edges of $G$, with $b_{ij}=\vartheta(v_i, e_j)$ (here $v_i\in V(G)$, $e_j \in E(G))$. Usually, when only $\Gamma$ is given, then we use the arbitrary orientation. So each row of the incidence matrix corresponding to vertex $v_i$ contains $d_{v_i}$ non-zero entries, each equal to $+1$ or $-1$. On the other hand, each column of the incidence matrix corresponding to edge $e_j$ contains two non-zero entries, each equal to $+1$ or $-1$. Therefore, even in the case that multiple edges exist, we easily obtain
 \begin{equation}
 B_{\vartheta}B_{\vartheta}^T = D(G)-A(\Gamma_{\vartheta})= L(\Gamma_{\vartheta}),
 \end{equation}
 where $D(G)$ is the diagonal matrix of vertex degrees of $G$. It is easy to observe that $L(\Gamma_{\vartheta}$) is positive-semidefinite.\\
  \indent The suddivision signed graph $A(S(\Gamma_{\vartheta}))$ is the signed graph whose underlying graph is $S(G)$ with vertex set $V(G) \cup E(G)$ . It preserves the orientation ${\vartheta}$ and can be represented in the block form as follows.
\begin{equation*}
  A(S(\Gamma_{\vartheta}))=
  \begin{pmatrix}
  O_n & B_{\vartheta}\\
  B_{\vartheta}^T & O_m
  \end{pmatrix},
\end{equation*}
 where $O_r \in M_r(\mathbb{R})$. It is easy to see that the signature $\sigma$ of the subdivision signed graph  is defined by $\sigma(v_ie_j)=\vartheta_{ij}$. An example of subdivision graph of a signed graph is shown in Figure 6.\\
{\bf Remark 4.1} Any  orientation (random) $\vartheta$  to the edges of $\Gamma$ gives rise to the same matrices $A(\Gamma_{\vartheta}) =A(\Gamma)$ and $L(\Gamma_{\vartheta}) =L(\Gamma)$, while the matrix $A(S(\Gamma_{\vartheta})) $ does depend on $\vartheta$. Let $S$ be a $\pm 1$ diagonal matrix such that $B_{\vartheta}'=B_{\vartheta} S$. It can be easily seen that $ A(S(\Gamma_{\vartheta'}))=[I_n \dot{+} S] A(S(\Gamma_{\vartheta}))[I_n \dot{+} S]$, where  $\dot{+}$ denotes the direct sum of two matrices. From now on, the index $\vartheta $ will be not specified anymore. \\
 \begin{figure}
\centering
	\includegraphics[scale=.8]{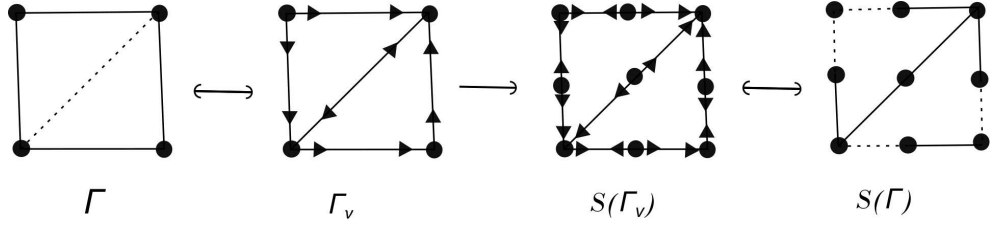}
	\caption{A signed graph and the corresponding signed subdivision  graph.}
	\label{Figure 6}
\end{figure}

\begin{lemma}\label{4.1}{\em \cite{bc}}
If $ {\mathit{B}}$ is the incident matrix of a connected signed graph $\Gamma=(G, \sigma)$ having $n$ vertices. Then
$$ rank(B)=\left\{\begin{array}{lr} n-1, &\mbox{if $\Gamma$ is balanced,}\\
n, &\mbox {if $\Gamma$ is unbalanced.}
\end{array}\right.$$
\end{lemma}
{\bf Operation 4.1} Let $\Gamma$ be a signed graph with vertex set $V(G)=\{v_1,v_2,\dots,v_n\}$ and edge set $E(G)=\{e_1,e_2,\dots,e_m\}$. Corresponding to each signed edge $e_i$, $1 \leq i \leq m$, in $\Gamma$, introduce a set $U_i^p$ of $p$ (positive integer) isolated vertices and make every vertex in $U_i^p$ adjacent to the vertices incident with $e_i$ , $i = 1, 2, \dots,m$ (in the same way as in subdivision signed graph $S(\Gamma)$) and remove edges of $\Gamma$ only. The resultant signed graph is denoted by $S_p(\Gamma)$.
 The number of vertices and edges of the signed graph $S_p(\Gamma)$ are $n+pm$ and $2pm$, respectively. If $p = 1$, then $S_p(\Gamma)$  coincides with the subdivision graph $S(\Gamma)$. Figure 7 illustrates the above operation.

\begin{theorem}\label{4.2}
Let $\Gamma$ be a signed graph with $n$ vertices and $m$ edges. Let $\mu_1\geq \mu_2\geq \dots \ge\mu_{n-1}>\mu_n\geq0$ be the Laplacian eigenvalues of the signed graph $\Gamma$. Then the adjacency spectrum of $S_p(\Gamma)$ is
$$ Spec(S_p(\Gamma))=\left\{\begin{array}{lr} \{0^{(pm-n+2)},\pm\sqrt{p\mu_1}^{(1)},\pm\sqrt{p\mu_2}^{(1)},\dots,\pm\sqrt{p\mu_{n-1}}^{(1)}\} &\mbox{if $\Gamma$ is balanced,}\\
\{0^{(pm-n)},\pm\sqrt{p\mu_1}^{(1)},\pm\sqrt{p\mu_2}^{(1)},\dots,\pm\sqrt{p\mu_{n-1}}^{(1)},\pm\sqrt{p\mu_{n}}^{(1)}\}, &\mbox {if $\Gamma$ is unbalanced.}
\end{array}\right.$$
\end{theorem}
{\bf Proof.} Using the suitable labelling of the vertices of $S_p(\Gamma)$, the adjacency matrix of $S_p(\Gamma)$ can be written as
\begin{equation*}
A(S_p(\Gamma))=
\begin{pmatrix}
   O &  B & B & \dots &  B\\
   B^T &  O & O & \dots &  O\\
   B^T &  O & O & \dots &  O\\
   \vdots  &  \vdots & \ddots &  \vdots\\
   B^T &  O & O & \dots &  O\\
\end{pmatrix}.
\end{equation*}
Therefore, we have
\begin{equation*}
A(S_p(\Gamma))^2=
\begin{pmatrix}
   O &  B & B & \dots &  B\\
   B^T &  O & O & \dots &  O\\
   B^T &  O & O & \dots &  O\\
   \vdots  &  \vdots & \ddots &  \vdots\\
   B^T &  O & O & \dots &  O\\
\end{pmatrix}
\begin{pmatrix}
   O &  B & B & \dots &  B\\
   B^T &  O & O & \dots &  O\\
   B^T &  O & O & \dots &  O\\
   \vdots  &  \vdots & \ddots &  \vdots\\
   B^T &  O & O & \dots &  O\\
\end{pmatrix}
\end{equation*}
\begin{figure}
\centering
	\includegraphics[scale=.8]{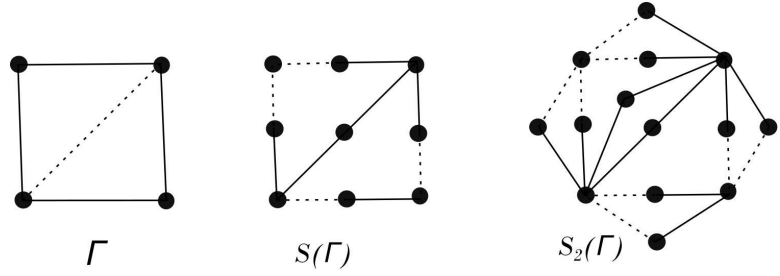}
	\caption{Signed graphs $\Gamma$, $S(\Gamma)$ and $S_2(\Gamma).$}
	\label{Figure 7}
\end{figure}
\begin{equation*}
  =\begin{pmatrix}
   pBB^T &  O & O & \dots &  O\\
   O &  B^TB & B^TB & \dots &  B^TB\\
   O &  B^TB & B^TB & \dots &  B^TB\\
   \vdots  &  \vdots & \ddots &  \vdots\\
   O &  B^TB & B^TB & \dots &  B^TB\\
\end{pmatrix}
\end{equation*}
\begin{equation*}
=\begin{pmatrix}
 pBB^T & O_{1\times p}\\
 O_{p\times 1} & J_{p\times p}\otimes B^TB\\
\end{pmatrix},
\end{equation*}
where  $J_{p\times p}\otimes B^TB $ denotes the Kronecker product of the matrices $ J_{p \times p} $ and $ B^TB $, and $ J_{p \times p} $ is a square matrix whose all entries are equal to $1$. Thus
\begin{equation*} Spec(A(S_p(\Gamma))^2)=  Spec(pBB^T) \cup Spec(J_{p\times p}\otimes B^TB).\end{equation*}
 As $ B^TB $ is a real symmetric matrix of order $m$, so all its eigenvalues are real. Let $x_1\geq x_2\geq \dots \geq x_m$ be the eigenvalues of the matrix $ B^TB $. Note that $rank(BB^T)=rank(B^TB)=rank(B)$. Therefore, by Lemma \ref{4.1}, we have
\begin{equation*} Spec(B^TB)=\left\{\begin{array}{lr} \{0^{(m-n+1)},x_1,x_2,\dots,x_{n-1}\} &\mbox{if $\Gamma$ is balanced,}\\
\{0^{(m-n)},x_1,x_2,\cdots,x_{n-1},x_{n}\} &\mbox {if $\Gamma$ is unbalanced}\end{array}\right. \end{equation*}
and
\begin{equation*} Spec(BB^T)=Spec(L(\Gamma))=\left\{\begin{array}{lr} \{0^{},\mu_1,\mu_2,\dots,\mu_{n-1}\} &\mbox{if $\Gamma$ is balanced,}\\
\{\mu_1,\mu_2,\dots,\mu_{n-1},\mu_{n}\} &\mbox {if $\Gamma$ is unbalanced,}\end{array}\right. \end{equation*}
where $x_n\neq 0$ and $\mu_n \neq 0$. As Spec($ J_{p \times p} $) is $ \{0^{p-1},p\}$, then by Lemma \ref{2.2}, we have
\begin{equation*} Spec(J_{p\times p}\otimes B^TB)=\left\{\begin{array}{lr} \{0^{(pm-n+1)},px_1,px_2,\cdots,px_{n-1}\} &\mbox{if $\Gamma$ is balanced,}\\
\{0^{(pm-n)},px_1,px_2,\dots,px_{n-1},px_{n}\} &\mbox {if $\Gamma$ is unbalanced.}\end{array}\right. \end{equation*}
We know that the underlying graph of a subdivision signed graph is always  bipartite. Therefore the underlying graph of $S_p(\Gamma)$ is always bipartite. Note that the eigenvalues of $B^TB$ are given by the eigenvalues of $BB^T$, together with $0$ of multiplicity $m-n$. Therefore, by Lemmas \ref{2.5} and  \ref{2.6}, we have
\begin{equation*} Spec(S_p(\Gamma))=\left\{\begin{array}{lr} \{0^{(pm-n+2)},\pm\sqrt{p\mu_1}^{(1)},\pm\sqrt{p\mu_2}^{(1)},\dots,\pm\sqrt{p\mu_{n-1}}^{(1)}\} &\mbox{if $\Gamma$ is balanced,}\\
\{0^{(pm-n)},\pm\sqrt{p\mu_1}^{(1)},\pm\sqrt{p\mu_2}^{(1)},\dots,\pm\sqrt{p\mu_{n-1}}^{(1)},\pm\sqrt{p\mu_{n}}^{(1)}\} &\mbox {if $\Gamma$ is unbalanced.}
\end{array}\right.\end{equation*}\qed
The following result can also be seen in [\cite{bc}, Theorem 2.2].
\begin{corollary}\label{4.3}
Let $\Gamma$ be a signed graph with $n$ vertices and $m$ edges. Let $\mu_1\geq \mu_2\geq \dots \ge\mu_{n-1}>\mu_n\geq0$ be the Laplacian eigenvalues of the signed graph $\Gamma$. Then the adjacency spectrum of $S(\Gamma)$ is
\begin{equation*} Spec(S(\Gamma))=\left\{\begin{array}{lr} \{0^{(m-n+2)},\pm\sqrt{\mu_1}^{(1)},\pm\sqrt{\mu_2}^{(1)},\dots,\pm\sqrt{\mu_{n-1}}^{(1)}\} &\mbox{if $\Gamma$ is balanced,}\\
\{0^{(m-n)},\pm\sqrt{\mu_1}^{(1)},\pm\sqrt{\mu_2}^{(1)},\dots,\pm\sqrt{\mu_{n-1}}^{(1)},\pm\sqrt{\mu_{n}}^{(1)}\} &\mbox {if $\Gamma$ is unbalanced.}
\end{array}\right.\end{equation*}
\end{corollary}
{\bf Operation 4.2} Let $\Gamma$ be a signed graph with vertex set $V(G)=\{v_1,v_2,\dots,v_n\}$ and edge set $E(G)=\{e_1,e_2,\dots,e_m\}$. Let $S(\Gamma)$ be a signed subdivision graph of a signed graph $\Gamma$ with vertex set $V(G)\cup E(G)$.  Corresponding to every vetrtex $v_i$, $1 \leq i \leq n$, in $S(\Gamma)$, introduce a set $V_i^p$ of $p$ (positive integer) isolated vertices and join each vertex of $V_i^p$ to the neighbors of $v_i$ with the same sign as the vertex $v_i$ in $S(\Gamma)$. Then in the resulting signed graph, corresponding to each vertex  $e_i$, ( $i = 1, 2, ...m$) introduce a set of $k$ isolated vertices $U_j^k$, $1\leq j\leq m$, where $k=p$ or $k=p-1$,  join each vertex in $U_j^k$ to the neighbors of $e_i$ with the same sign as the vertex $e_i$ in $S(\Gamma)$. The resultant signed graph is denoted by $S_p^k(\Gamma)$. The number of vertices and edges of the graph $S_p^k(\Gamma)$ are $(p+1)n+(k+1)m$ and $2(p+1)(k+1)m$, respectively. Figure 8 illustrates the above operation.
\begin{theorem}\label{4.4}
Let $\Gamma$ be a signed graph with $n$ vertices and $m$ edges. Let $\mu_1\geq \mu_2\geq \dots \ge\mu_{n-1}>\mu_n\geq0$ be the Laplacian eigenvalues of the signed graph $\Gamma$. Then the adjacency spectrum of $S_p^k(\Gamma)$ is\\
$ Spec(S_p^k(\Gamma))= \{0^{((p-1)n+(k+1)m+2)},\pm\sqrt{(p+1)(k+1)\mu_1}^{(1)},\dots,\pm\sqrt{(p+1)(k+1)\mu_{n-1}}^{(1)}\}$ \\if $\Gamma$ is balanced, and\\
$Spec(S_p^k(\Gamma))=\{0^{((p-1)n+(k+1)m)},\pm\sqrt{(p+1)(k+1)\mu_1}^{(1)},\dots,\pm\sqrt{(p+1)(k+1)\mu_{n-1}}^{(1)},\\ ~~~~~~~~~~~~~~~~ \pm\sqrt{(p+1)(k+1)\mu_{n}}^{(1)}\}$  \\if $\Gamma$ is unbalanced.
\end{theorem}
{\bf Proof.} By the suitable labelling of the vertices of $S_p^k(\Gamma)$, the adjacency matrix  $S_p^k(\Gamma)$ can be written as
\begin{equation*}
A(S_p^k(\Gamma))=
\begin{pmatrix}
   O &  B & O & \dots & O &  B\\
   B^T &  O & B^T & \dots & B^T &  O\\
   O &  B & O & \dots & O &  B\\
   \vdots  &  \vdots & \ddots &  \vdots\\
   O & B & O & \dots &O & B\\
   B^T &  O & B^T  & \dots & B^T &  O\\
\end{pmatrix},
\end{equation*}
when $k=p$. If $k=p-1$, then we have
\begin{equation*}
A(S_p^k(\Gamma))=
\begin{pmatrix}
   O &  B & O & \dots & B &  O\\
   B^T &  O & B^T & \dots & O &  B^T\\
   O &  B & O & \dots & B &  O\\
   \vdots  &  \vdots & \ddots &  \vdots\\
   B^T & O & B^T & \dots &O & B^T\\
   O &  B & O  & \dots & B &  O\\
\end{pmatrix}.
\end{equation*}
To prove the result, the following two cases arise.\\
{\bf Case 1.} Let $\Gamma$ be a balanced signed graph with $n$ vertices and $m$ edges. Let $X\in M_{n\times 1}(\mathbb{R})$ and $Y\in M_{m\times 1}(\mathbb{R})$ be two non-zero column vectors. Let $ Z= \begin{pmatrix}
X\\Y
\end{pmatrix} \in M_{(n+m)\times 1}(\mathbb{R}) $ be the eigenvector corresponding to the non-zero eigenvalue $\lambda_i$, $1\leq i\leq 2n-2$, of $S(\Gamma)$. Then $A(S(\Gamma))Z= \lambda_i Z$ implies that $BY=\lambda_i X$ and $B^TX=\lambda_i Y $. To find the eigenvalues of $S_p^k(\Gamma)$, consider the following two subcases.\\
\begin{figure}
\centering
	\includegraphics[scale=.8]{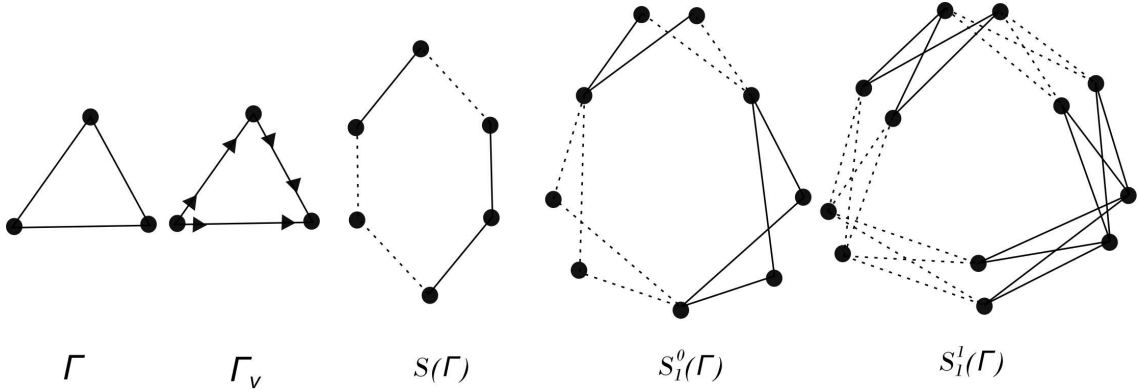}
	\caption{Signed graphs $\Gamma$, $\Gamma_{\vartheta}$, $S(\Gamma)$, $S_1^0(\Gamma)$  and $S_1^1(\Gamma).$}
	\label{Figure 8}
\end{figure}
{\bf Subcase 1.1.} If $k=p$, then let $U=\begin{pmatrix}
X\\Y\\ \vdots\\X\\Y
\end{pmatrix} \in M_{((p+1)n+(k+1)m)  \times 1}(\mathbb{R})$ be a non-zero column vector such that
\begin{equation*}
 A(S_p^k(\Gamma))U = \begin{pmatrix}
 O &  B & O & \dots & O &  B\\
   B^T &  O & B^T & \dots & B^T &  O\\
   O &  B & O & \dots & O &  B\\
   \vdots  &  \vdots & \ddots &  \vdots\\
   O & B & O & \dots &O & B\\
   B^T &  O & B^T  & \dots & B^T &  O\\
\end{pmatrix}
\begin{pmatrix}
X\\Y\\ \vdots\\X\\Y
\end{pmatrix}
=\begin{pmatrix}
(p+1)\lambda_i X\\ (p+1)\lambda_i Y \\ \vdots\\ (p+1)\lambda_iX\\ (p+1)\lambda_iY
\end{pmatrix}
\end{equation*}
\begin{equation*}=(p+1)\lambda_i U.\end{equation*}
Therefore $(p+1)\lambda_i$ is an eigenvalue of 	$S_p^k(\Gamma)$ corresponding to an eigenvector $U$. Thus the result follows by Corollary \ref{4.3}.\\
{\bf Subcase 1.2.} If $k=p-1$, then let $U=\begin{pmatrix}
\sqrt{k+1}X\\\sqrt{p+1}Y\\ \vdots\\\sqrt{k+1}X\\ \sqrt{p+1}Y\\\sqrt{k+1}X
\end{pmatrix} \in M_{((p+1)n+(k+1)m ) \times 1}(\mathbb{R})$ be a non-zero column vector such that
\begin{equation*}
 A(S_p^k(\Gamma))U = \begin{pmatrix}
  O &  B & O & \dots & B &  O\\
   B^T &  O & B^T & \dots & O &  B^T\\
   O &  B & O & \dots & B &  O\\
   \vdots  &  \vdots & \ddots &  \vdots\\
   B^T & O & B^T & \dots &O & B^T\\
   O &  B & O  & \dots & B &  O\\
\end{pmatrix}
\begin{pmatrix}
\sqrt{k+1}X\\\sqrt{p+1}Y\\ \vdots\\\sqrt{k+1}X\\ \sqrt{p+1}Y\\\sqrt{k+1}X
\end{pmatrix}
\end{equation*}
\begin{equation*}
=\begin{pmatrix}
(k+1)\lambda_i\sqrt{p+1}X\\ (p+1)\lambda_i\sqrt{k+1}Y\\ \vdots\\ (k+1)\lambda_i\sqrt{p+1}X\\ (p+1)\lambda_i \sqrt{k+1}Y\\(k+1) \lambda_i\sqrt{p+1}X
\end{pmatrix}
\end{equation*}
\begin{equation*}=\sqrt{(k+1)(p+1)}\lambda_i U.\end{equation*}
Therefore $\sqrt{(k+1)(p+1)}\lambda_i$ is an eigenvalue of 	$S_p^k(\Gamma)$ corresponding to an eigenvector $U$. Therefore the result follows by Corollary \ref{4.3}.\\
{\bf Case 2.}  When $\Gamma$ is an unbalanced signed graph with $n$ vertices and $m$ edges, the proof is similar to that  of Case $1$.\qed

Various constructions for cospectral non-isomorphic   regular graphs, cospectral non-isomorphic  Laplacian graphs and cospectral non-isomorphic  signless Laplacian graphs can be seen in {\cite{1,2,3,4,5,6}}. The following results shows that these constructions including the constructions obtained in the last section can be utilized to obtain infinite families of cospectral non-isomorphic  signed graphs.
\begin{corollary}\label{4.5}
Let $\Gamma_1$ and $\Gamma_2$ be two non-isomorphic signed graphs which are Laplacian cospectral. Then\\
$(i)$ the signed graphs $S_p(\Gamma_1)$ and $S_p(\Gamma_2)$ are cospectral and non-isomorphic,\\
$(ii)$ the signed graphs $S_p^k(\Gamma_1)$ and $S_p^k(\Gamma_2)$ are cospectral and non-isomorphic.
\end{corollary}
{\bf Proof.} Let $\Gamma_1$ and $\Gamma_2$ be two non-isomorphic signed graphs. Then, clearly $S_p(\Gamma_1)$ and $S_p(\Gamma_2)$ are non-isomorphic signed graphs and $S_p^k(\Gamma_1)$ and $S_p^k(\Gamma_2)$ are non-isomorphic signed graphs. Hence the result follows by Theorems \ref{4.2} and \ref{4.4}.\\
\noindent{\bf Example 4.1}
Consider the two non-isomorphic signed graphs $\Gamma_1$ and $\Gamma_2$, which are Laplacian cospectral, as shown in Figure $9$. Their Laplacian spectrum is respectively given by
$Spec_L(\Gamma_1)=\{0,2,3^{(2)},3+\sqrt{
5}, 3-\sqrt{
5}\}$ and $Spec_L(\Gamma_2)=\{0,2,3^{(2)},3+\sqrt{
5}, 3-\sqrt{
5}\}$.
It is easy to see that $S_2(\Gamma_1)$ and $S_2(\Gamma_2)$ are non-isomorphic signed graphs which are cospectral as their adjacency spectrum are respectively given by $Spec(S_2(\Gamma_1))=\{0^{(10)},\pm 2,\pm  \sqrt{6}^{(2)},\pm(\sqrt{6+\sqrt{
20})}, \pm(\sqrt{6-\sqrt{
20})}\}$ and $Spec(S_2(\Gamma_2))=\{0^{(10)},\pm 2,\pm  \sqrt{6}^{(2)},\pm(\sqrt{6+\sqrt{
20})}, \pm(\sqrt{6-\sqrt{
20})}\}$.
\begin{figure}
\centering
	\includegraphics[scale=.8]{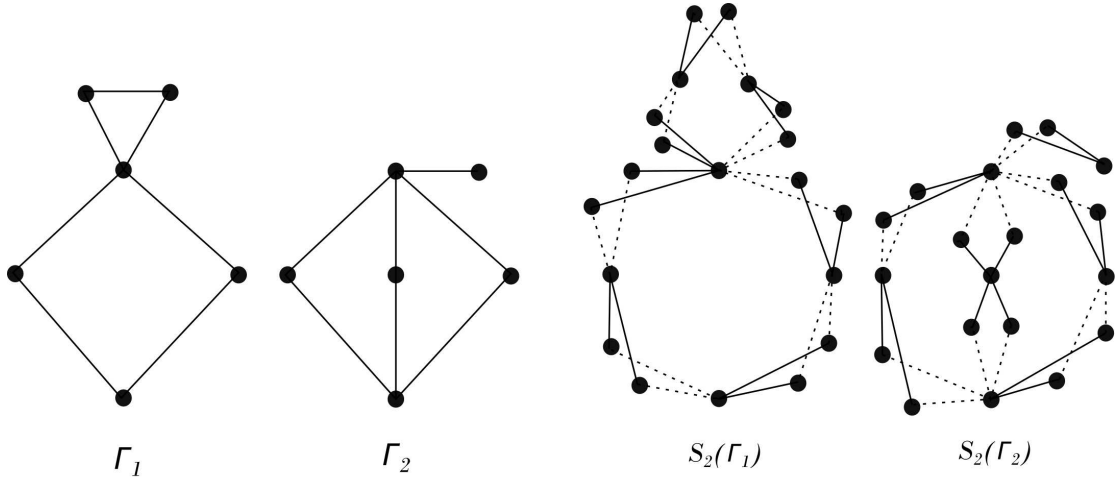}
	\caption{Cospectral signed graphs $S_2(\Gamma_1)$ and $S_2(\Gamma_2)$.}
	\label{Figure 9}
\end{figure}

\begin{corollary}\label{4.6}
Let $\Gamma_1$ and $\Gamma_2$ be two cospectral non-isomorphic  $r-$regular signed graphs. Then\\
$(i)$ the signed graphs $S_p(\Gamma_1)$ and $S_p(\Gamma_2)$ are cospectral and non-isomorphic,\\
$(ii)$ The signed graphs $S_p^k(\Gamma_1)$ and $S_p^k(\Gamma_2)$ are cospectral and non-isomorphic.
\end{corollary}
{\bf Proof.} Let $\Gamma_1$ and $\Gamma_2$ be two non-isomorphic cospectral regular signed graphs. Then $L(\Gamma_1)=D(\Gamma_1)-A(\Gamma_1)$ and $L(\Gamma_2)=D(\Gamma_2)-A(\Gamma_2)$ are cospectral. Hence the result follows by Corollary \ref{4.5}.\qed

\begin{corollary}\label{4.7}
Let $\Gamma$ be a signed graph whose all Laplacian eigenvalues are perfect squares. Then\\
$(i)$ the signed graph $S_p(\Gamma)$ is integral, if p is a perfect square,\\
$(ii)$ the signed graph $S_p^k(\Gamma)$ is integral, if (p+1)(k+1) is a perfect square.
\end{corollary}

\noindent{\bf Example 4.2} Let $K_n$ be a balanced complete signed graph on $n$ vertices, where $n=t^2$, $t\geq 2$ is a positive integer. Then\\ $(i)$ the signed graph $S_p(K_n)$ is integral, if $p$ is a perfect square,\\
$(ii)$  the signed graph $S_p^k(K_n)$ is integral, if $(p+1)(k+1)$ is a perfect square.\\
The following result is the graceful implication of Lemma \ref{2.7} and Corollaries \ref{4.5} and \ref{4.6}.
\begin{theorem}\label{4.8}  For infinitely many $n$, there exists a family of $2^k$ pairwise
cospectral nonisomorphic  signed graphs on $n$
vertices, where $k > \frac{n}{•(2 log_2(n))}$.
\end{theorem}

The following result directly follows from Theorems \ref{4.2} and \ref{4.4}.
\begin{theorem}\label{4.9}
Let $\Gamma$ be a signed graph with $n$ vertices and $m$ edges. Then    \\$(i)$ $\mathcal{E}(S_p(\Gamma))=\sqrt{p}\mathcal{E}(S(\Gamma)),$\\
$(ii)$ $\mathcal{E}(S_p^k(\Gamma))=\sqrt{(p+1)(k+1)}\mathcal{E}(S(\Gamma)).$
\end{theorem}
\begin{theorem}\label{4.10}
Let $\Gamma$ be an unbalanced unicyclic signed graph with at least one edge and having  Laplacian eigenvalues $\mu_1 \geq \mu_2\geq \cdots \geq \mu_n>0$. Then $S(\Gamma)\times K_2$ and $S(\Gamma)\otimes K_2$ are  noncospectral and equienergetic if and only if $\mu_n \ge 1$.
\end{theorem}
{\bf Proof.} Let $\Gamma$ be an unbalanced unicyclic signed graph. Then, by Theorem \ref{4.2}, we have
\begin{equation*} Spec(S(\Gamma))=\{\pm\sqrt{\mu_1}^{(1)},\pm\sqrt{\mu_2}^{(1)},\dots,\pm\sqrt{\mu_{n-1}}^{(1)},\pm\sqrt{\mu_{n}}^{(1)}\}. \end{equation*}
First, assume that  $\mu_n \ge 1$. This implies that $|\sqrt{\mu_j}| \ge 1$, for all $j=1,2,\dots,n$.  Also,
\begin{equation*}E(S(\Gamma)\times K_2)=2\sum\limits_{j=1}^{n}(|\sqrt{\mu_j}+1|+|\sqrt{\mu_j}-1|).\end{equation*}
\indent As $|\sqrt{\mu_j}|\ge 1$, for all $j=1,2,\dots,n$, we have\\
\begin{align*}
E(S(\Gamma)\times K_2)&=2\sum\limits_{j=1}^{n}(|\sqrt{\mu_j}|+1+|\sqrt{\mu_j}|-1)\\&=2E(S(\Gamma))\\&=E(S(\Gamma))E(K_2)=E(S(\Gamma)\otimes K_2).
\end{align*}
\indent Note that $\sqrt{\mu_1}+1\in Spec(S(\Gamma)\times K_2)$ but $\sqrt{\mu_1}+1\notin Spec(S(\Gamma)\otimes K_2)$. Therefore $S(\Gamma)\times K_2$ and $S(\Gamma)\otimes K_2$ are noncospectral.  The converse is similar to that of the converse in Lemma \ref{2.4}.\qed

\noindent{\bf Example 4.3.} Let $C_3^- =(C_3,-)$ be an unbalanced unicyclic signed graph on $3$ vertices. Its Laplacian spectrum is given by $Spec_L(C_3^-)=\{4,1,1\}$. Therefore $C_3^-$ meets the requirement of Theorem $4.10$. Hence  $S(C_3^-)\times K_2$ and $S(C_3^-)\otimes K_2$ are  noncospectral and equienergetic.\\
\begin{figure}
\centering
	\includegraphics[scale=.8]{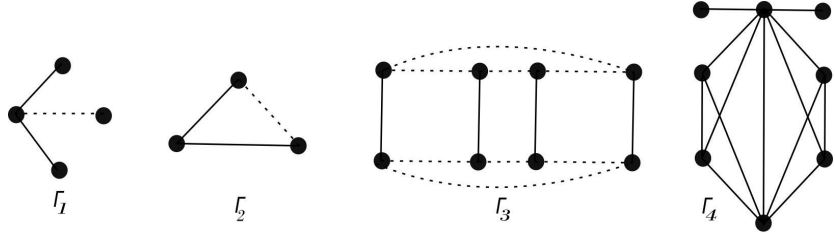}
	\caption{Signed graphs $\Gamma_1$,  $\Gamma_2$, $\Gamma_3$ and $\Gamma_4$.}
	\label{Figure 10}
\end{figure}

The following corollary directly follows from Theorems \ref{4.2}, \ref{4.4} and \ref{4.9}.
\begin{corollary}\label{4.11}
Let $\Gamma_1$ and $\Gamma_2$ be two signed graphs whose signed subdivision graphs are noncospectral and equienergetic. Then \\
$(i)$ the signed graphs $S_p(\Gamma_1)$ and $S_p(\Gamma_2)$are noncospectral and equienergetic,\\
$(ii)$ the signed graphs $S_p^k(\Gamma_1)$ and $S_p^k(\Gamma_2)$are noncospectral and equienergetic.
\end{corollary}
\noindent{\bf Example 4.4.} Consider the signed graphs $\Gamma_1$, $\Gamma_2$, $\Gamma_3$ and $\Gamma_4$ as shown in Figure 10. The adjacency spectrum of their signed subdivision graphs is respectively given by $ Spec(S(\Gamma_1))=\{\pm 2, \pm 1^{(2)},0\} $, $ Spec(S(\Gamma_2))=\{\pm 2,\pm 1^{(2)}\} $, $ Spec(S(\Gamma_3))=\{\pm 2^{(3)},\pm \sqrt{2}^{(3)},\pm \sqrt{6}, 0^{(6)}\} $ and\\ $ Spec(S(\Gamma_4))=\{\pm 1^{(2)},\pm 2^{(2)},\pm \sqrt{2}^{},\pm 2\sqrt{2},\pm \sqrt{6}, 0^{(7)}\}. $  Clearly, the signed graphs $S(\Gamma_1)$ and $S(\Gamma_2)$ are noncospectral and equienergetic. Similarly, the signed graphs $S(\Gamma_3)$ and $S(\Gamma_4)$ are noncospectral and equienergetic. Thus, by  Corollary \ref{4.11}, we have\\
$(i)$ the signed graphs $S_p(\Gamma_1)$ and $S_p(\Gamma_2)$are noncospectral and equienergetic,\\
$(ii)$ the signed graphs $S_p^k(\Gamma_1)$ and $S_p^k(\Gamma_2)$are noncospectral and equienergetic,\\
$(iii)$ the signed graphs $S_p(\Gamma_3)$ and $S_p(\Gamma_4)$are noncospectral and equienergetic,\\
$(iv)$ the signed graphs $S_p^k(\Gamma_3)$ and $S_p^k(\Gamma_4)$are noncospectral and equienergetic.\\

\noindent{\bf Acknowledgements.}   This research is supported by SERB-DST research project number CRG/2020/000109. The research of Tahir Shamsher is supported by SRF financial assistance by Council of Scientific and Industrial Research (CSIR), New Delhi, India.

\noindent{\bf Data availibility} Data sharing is not applicable to this article as no datasets were generated or analyzed
during the current study.


\begin{thebibliography}{0}

\bibitem{bc} F. Belardo and S. K.Simic, On the Laplacian coefficients of signed graphs, Linear Algebra Appl. {\bf 475} (2015) 94-113.
\bibitem{bcs}  F. Belardo, M. Brunetti, M. Cavaleri and A. Donno, Constructing cospectral signed graphs, Linear Multilinear Algebra (2019) https://doi.org/10.1080/03081087.2019.1694483.
\bibitem{openproblem} F. Belardo, S. M.Cioaba, J. Koolen and J. Wang, Open problems in the spectral theory of
signed graphs, The Art of Discrete and Applied Mathematics 1 (2018) P2.10 https://doi.org/10.26493/2590-9770.1286.d7b.
\bibitem{mpe} M. A. Bhat and S. Pirzada, On equienergetic signed graphs, Discrete. Appl. Math. {\bf 189} (2015) 1-7.
\bibitem{mpt} M. A. Bhat and S. Pirzada, Unicyclic signed graphs with minimal energy, Discrete. Appl. Math. {\bf 226} (2017) 32-39.
\bibitem{1} Z. L. Blazsik, J. Cummings, and  W. H. Haemers, Cospectral regular graphs with and
without a perfect matching, Discrete Math. {\bf 338}(2015) 199-201.
\bibitem{ds} D. Cvetkovic, S. Simic, and P. Rowlinson, An Introduction to the Theory of Graph Spectra, Cambridge University Press, 2009.
\bibitem{ghz} K. A. Germina, S. Hameed  and T. Zaslavsky, On products and line graphs of signed graphs, their eigenvalues and energy, Linear Algebra Appl. {\bf 435} (2010) 2432-2450.

\bibitem{2} A. Dehghan and A. H. Banihashemi,  Cospectral bipartite graphs with the same
degree sequences but with different number of large cycles, Graphs Combin. {\bf 35} (2019) 1673-1693.
\bibitem{3} S. Dutta,  Constructing non-isomorphic signless Laplacian cospectral graphs, Discrete Math. {\bf 343} (2020)111783.
\bibitem{4} C. D. Godsil and B. D. McKay, Constructing cospectral graphs, Aequationes Math.{\bf 25} (1982) 257-268.
\bibitem{hj} R. A. Horn and C. R. Johnson, Matrix Analysis,
Cambridge University Press, 23-Feb-1990.
\bibitem{h} F. Harary, On the notion of balanced in a signed graph, Michigan Math. J. {\bf 2}
(1953) 143-146.
\bibitem{6} M. Haythorpe and  A. Newcombe, Constructing families of cospectral regular graphs,
Combin. Probab. Comput. {\bf 29} (2020) 664-671.
\bibitem{5}  W. H. Haemers and E. Spence, Enumeration of cospectral graphs, European J. Combin. {\bf{ 25}} (2004)199-211.
\bibitem{y} Y. Hou, J. Li and Y. Pan, On the Laplacian eigenvalues of signed graphs. Linear Multilinear Algebra {\bf 51} (2003) 21-30.
\bibitem{r} R. Merris, Large Families of Laplacian isospectral graphs, Linear Multilinear Algebra,  {\bf 43} (1997) 201 -205.
\bibitem{Z} T. Zaslavsky, A mathematical bibliography of signed and gain graphs and allied areas, Electron. J. Combin., Dyn. Surv. (1999) DS8, http://www3.combinatorics.org/Surveys/ds8.pdf.


\end{thebibliography}
\end{document}